\documentclass{anstrans}
\title{A Linear, Exponential-Discontinuous Scheme for Discrete-Ordinates Calculations in Slab Geometry}
\author{Jeremy A. Roberts$^{*}$}

\institute{
$^{*}$Department of Mechanical and Nuclear Engineering, Kansas State University,
Manhattan, KS, jaroberts@ksu.edu
}

\usepackage{graphicx} 
\usepackage{booktabs} 
\usepackage{microtype} 
\usepackage{braket}

\begin{document}
\section{Introduction}

In the development of spatial discretization schemes for the discrete-ordinates (S$_N$) equations, one often seeks some combination of accuracy, efficiency, and the preservation of certain physical characteristics like the positivity of the angular flux solution.  The problem of accuracy is reasonably straightforward to address for optically-thin cells, i.e., when the size of the cell is much smaller than a mean-free path.  However, for optically-thick cells, accuracy can suffer greatly if the underlying solution is diffusive\cite{larsen1987asn} and the scheme does not have the thick-diffusion limit.  Moreover, most schemes can produce negative fluxes, with the step-characteristic schemes being a well-known exception\cite{lathrop1969sdt}.  The problem of negative fluxes is often exacerbated in the presence of optically-thick cells, which may be required for the analysis of deep-penetration problems subject to limited computational resources (e.g., when using deterministic models to generate variance-reduction parameters for Monte Carlo models).

Presented here is a preliminary study of a strictly linear, discontinuous-Petrov-Galerkin scheme for the discrete-ordinates method in slab geometry.  By ``linear'', we mean the  discretization does not depend on the solution itself as is the case in classical ``fix-up'' schemes\cite{lewis1984cmn} and other nonlinear schemes that have been explored to maintain positive solutions with improved accuracy (e.g., \cite{walters1996asp, mathews1994ecs, maginot2012nnm}).  By discontinuous, we mean the angular flux $\psi$ and scalar flux $\phi$ are piecewise continuous functions that may exhibit discontinuities at cell boundaries.  Finally, by ``Petrov-Galerkin,'' we mean a finite-element scheme in which the ``trial'' and ``test'' functions differ.  In particular, we find that a trial basis consisting of a constant and exponential function that exactly represents the step-characteristic solution with a constant and linear test basis produces a scheme (1) with slightly better local errors than the linear-discontinuous (LD) scheme (for thin cells), (2) accuracy that approaches the linear-characteristic (LC) scheme (when the LC solution is positive), and (3) is positive as long as the first two source Legendre moments satisfy $|s_1| < 3 s_0$.

\section{Theory}

The discrete-ordinates equations with isotropic scattering and sources in slab geometry are
\begin{equation}
  \mu_m \frac{d\psi_m}{dx} + \Sigma_t(x) \psi_n(x) = \frac{1}{2}\left (\Sigma_s(x) \phi(x) + s(x) \right ) \, ,
\end{equation}
and
\begin{equation}
  \phi(x) = \int^1_{-1} \psi_m(x)d\mu \approx \sum^M_{m=1} w_m \psi_m(x) \, ,
\end{equation}
where $\psi_m(x) = \psi(x, \mu_m)$ is the angular flux for a discrete direction with cosine $\mu_m$ relative to the $x$ axis, and the notation is otherwise standard.

To simplify the notation a bit, we drop the $m$ subscript for $\mu$, use $s(x)$ for the entire right-hand side, and limit the spatial domain to a single cell of width $\Delta$, or
\begin{equation}
\mu \frac{d\psi}{dx} + \Sigma \psi(x) = s(x) \, , \quad x \in [0, \Delta] \, ,
\label{eq:simple_te}
\end{equation}
with $\psi(0)$ or $\psi(\Delta)$ given for $\mu > 0$ and $\mu < 0$; the presentation will assume $\mu > 0$, but generalization to $\mu < 0$ is straightforward. As part of a larger model, one can consider this cell to be the $k$th cell of several representing a larger domain.  In such a case, the given boundary fluxes represent fluxes emerging from neighboring cells $k-1$ and $k+1$ or a global boundary.

\subsection{The DG Framework: A Quick Review}

Our basic approach for solving Eq. (\ref{eq:simple_te}) numerically is to assume that
\begin{equation}
  \psi(x) \approx \tilde{\psi} \equiv \sum^N_{i=0} \psi_i p^{\psi}_i(x) \, ,
\end{equation}
and
\begin{equation}
  s(x) \approx \tilde{s} \equiv \sum^N_{i=0} s_i p^s_i(x) \, ,
\end{equation}
where the bases $\mathcal{P}_{\psi} = [p^{\psi}_0(x), \ldots, p^{\psi}_N(x)]$ and $\mathcal{P}_{s} = [p^{s}_0(x), \ldots, p^{s}_N(x)]$ are often (but not required to be) identical sets of polynomials orthogonal over $x \in [0, \Delta]$.  In the finite-element world, these functions (well, specifically those in $\mathcal{P}_{\psi}$) are called ``trial functions.''  (For neutron transport, the terminology is a bit confused because we solve for $\psi(x)$ in a single cell, but that solution contributes to $\phi(x)$, and, hence, $s(x)$!)

Now, define a ``test'' space spanned by $\mathcal{P} = [p_0(x), \ldots, p_N(x)]$.  Then, substitution of $\tilde{\psi}(x)$ into Eq.~ (\ref{eq:simple_te}), multiplication of the result by the
``test function'' $p_i(x)$, and integration of that product over the domain leads to
\begin{equation}
\begin{split}
\mu \int^{\Delta}_0 p_i(x) \frac{d\tilde{\psi}}{dx} dx
    &+\Sigma \int^{\Delta}_0 p_i(x) \tilde{\psi}(x) dx \\
    &=  \int^{\Delta}_0 p_i(x) \tilde{s}(x) dx \, ,  \quad i = 0, 1, \ldots N \, ,
\end{split}
\end{equation}
or, following the application of integration by parts of the first term, we have
\begin{equation}
\begin{split}
\mu p_i(\Delta) \tilde{\psi}(\Delta)
  &- \mu p_i(0) \hat{\psi}(0)
  - \mu \int^{\Delta}_0 \frac{dp_i}{dx} \tilde{\psi}(x) dx \\
  &+ \Sigma \int^{\Delta}_0 p_i(x) \tilde{\psi}(x) dx \\
    &= \int^{\Delta}_0 p_i(x) \tilde{s}(x) dx \, ,  \quad i = 0, 1, \ldots N \, ,
\end{split}
\label{eq:projected_general}
\end{equation}
which represents $N+1$ equations for the flux moments $\psi_i$.  Certain choices for the bases (e.g., polynomials) greatly simplify these equations.

Note the use of $\hat{\psi}$ at the left boundary: $\hat{\psi}(0)$ is the (effective) incident angular flux $\hat{\psi}(0)$ at the left edge of cell $k$ given the edge flux on either side of the edge, i.e., $\tilde{\psi}_{k-1}(0^-)$ and $\tilde{\psi}_{k}(0^+)$ for $\mu > 0$.  Importantly,  {\it discontinuous-Galerkin methods} as introduced by Reed and Hill\cite{reed1973triangular} allow us to choose $\hat{\psi}(0)$ for $\mu > 0$ such that the inflow of neutrons is not represented exactly by the approximation $\tilde{\psi}(x)$, $x \in [0, \Delta]$.  This mismatch is the ``discontinuous'' in DG.
It is common to define effective (more often called ``numerical'') fluxes using the generic form
\begin{equation}
 \hat{\psi}(0) = f(\tilde{\psi}_{k-1}(0^-), \tilde{\psi}_{k}(0^+)) \, ,
\end{equation}
but we'll stick with the {\it upwind approximation}, i.e., we set
\begin{equation}
 \hat{\psi}(0) = \tilde{\psi}_{k-1}(0^-) \,
\end{equation}
for $\mu > 0$.

Perhaps the best-known discontinuous-Galerkin method is the linear-discontinuous (LD) approximation in which the flux and source are represented by a constant and linear term (although an equivalent representation using two linear ``cardinal'' functions that are unity at one edge and zero at the other is sometimes more convenient\cite{adams2001discontinuous}).
To illustrate the mechanics, we'll use the scaled and shifted Legendre polynomials, the first two of which are
\begin{equation}
 p_{0}(x) = 1 \quad \text{and} \quad p_{1}(x) = \frac{2x - \Delta}{\Delta} \, .
\end{equation}
Application of Eq.~(\ref{eq:projected_general}) for $N=1$ with these  basis functions leads to
\begin{equation}
\left [ \begin{matrix}
  \Delta \Sigma + \mu & \mu \\
  -\mu                & \frac{\Delta \Sigma}{3} + \mu \\
\end{matrix} \right ]
\left [ \begin{matrix}
 \psi_0 \\
 \psi_1 \\
\end{matrix} \right ] =
\left [ \begin{matrix}
 \Delta s_0 + \mu \hat{\psi} \\
  \frac{\Delta s_1 }{3} -\mu \hat{\psi} \\
\end{matrix} \right ]  \, .
\end{equation}
The corresponding solution is $\psi_{LD}(x) = \psi_0 p_0(x) + \psi_1 p_1(x)$.  For comparison later on, we note that, for a strictly linear source, the local errors in the incident, outgoing, and average angular flux are
\begin{equation}
  \epsilon_{in} = - \frac{s_{1} \tau}{3 \Sigma \mu} + \tau^{2} \left(- \frac{\psi_{up}}{6 \mu^{2}} + \frac{s_{0}}{6 \Sigma \mu^{2}} + \frac{s_{1}}{18 \Sigma \mu^{2}}\right) + O\left(\tau^{3}\right) \, ,
  \label{eq:ld_error_in}
\end{equation}
\begin{equation}
  \epsilon_{out} = - \frac{s_{1} \tau^{3}}{36 \Sigma \mu^{3}} + \tau^{4} \left(- \frac{\psi_{up}}{72 \mu^{4}} + \frac{s_{0}}{72 \Sigma \mu^{4}} + \frac{23 s_{1}}{1080 \Sigma \mu^{4}}\right) + O\left(\tau^{5}\right),
  \label{eq:ld_error_out}
\end{equation}
and
\begin{equation}
  \epsilon_{avg} = \frac{s_{1} \tau^{2}}{36 \Sigma \mu^{2}} + \tau^{3} \left(\frac{\psi_{up}}{72 \mu^{3}} - \frac{s_{0}}{72 \Sigma \mu^{3}} - \frac{23 s_{1}}{1080 \Sigma \mu^{3}}\right) + O\left(\tau^{4}\right) \, ,
  \label{eq:ld_error_avg}
\end{equation}
where $\tau = \Sigma\Delta$ is the optical cell width.

\subsection{Mixed Bases: A DPG Framework}

In the LD and step methods, the flux and source are represented as Legendre polynomials.  Now, suppose that the flux were represented in one basis and the source in another.  The basis functions are not restricted. It helps to use the constant function for the first basis function to preserve averages, but this is not enforced.   Further, it helps to use basis functions $p_i(x)$ that are mutually orthogonal over the cell, i.e.,
\begin{equation}
  \int^{\Delta}_0 p_i(x) p_j(x) dx = 0 \, , \quad  i \neq j \, .
\end{equation}
We will assume such orthogonality as a matter of convenience.

For brevity, let us define the inner product
\begin{equation}
  \braket{f_i| g_j} = \int^{\Delta}_{0} f_i(x) g_j(x) dx  \, ,
\end{equation}
so that Eq.~(\ref{eq:projected_general}) becomes
\begin{equation}
\mu p_i(\Delta) \tilde{\psi}(\Delta)
  - \mu p_i(0) \hat{\psi}(0)
  - \mu \braket{\frac{dp_i}{dx}|\tilde{\psi}}
  + \Sigma \braket{p_i|\tilde{\psi}}
   =  \braket{p_i|\tilde{s}} \, .
\label{eq:projected_general_braket}
\end{equation}
Notice that the test functions $p_i(x)$ are not assumed to be equal to the flux or source basis functions.  If $p_i(x) = p^s_i(x), \forall i$, then the right-hand side simplifies to $s_i$.  If $p_i(x) = p^{\psi}(x)$, then the last term on the left simplifies to $\Sigma \psi_i$.  When the test functions differ from the trial functions, the resulting scheme is more accurately called a discontinuous-Petrov-Galerkin (DPG) method.  There is an entire class of methods known as DPG methods in the literature (see, e.g., \cite{demkowicz2010class}); however, beyond the use of mixed bases first suggested by Petrov for Galerkin methods\cite{petrov1940application}, there appears to be no obvious connection.

If we choose $\mathcal{P}_{\psi} \neq \mathcal{P}_{s}$, we must project the resulting flux (in its basis) onto the source basis in order to define scalar flux-moment contributions.  Let $\boldsymbol{\psi} = [\psi_0, \psi_1, \ldots]^T$ be the flux moments in the flux basis; we seek the corresponding moments $\boldsymbol{\varphi} = [\varphi_0, \varphi_1, \ldots]^T$ in the source basis.
For the $i$th flux basis function $p^{\psi}_i(x)$, we can write
\begin{equation}
  p^{\psi}_i(x) \approx \sum_{j=0}^{N} \frac{ \braket{ p^s_j | p^{\psi}_i }}{\braket{p^s_j | p^s_j }} p^s_j(x) \, .
\end{equation}
Hence,
\begin{equation}
\begin{split}
  \tilde{\psi} &\approx \sum^N_{i=0} \psi_i \sum_{j=0}^{N} \frac{\braket{p^s_j| p^{\psi}_i}}{\braket{p^s_j|p^s_j}} p^s_j(x) \\
   &=  \sum_{j=0}^{N}   \left ( \sum^N_{i=0} \psi_i  \frac{\braket{p^s_j| p^{\psi}_i}}{\braket{p^s_j|p^s_j}} \right ) p^s_j(x) \\
   &=  \sum_{j=0}^{N} \varphi_j p^s_j(x) \, ,
\end{split}
\end{equation}
It follows that
\begin{equation}
  \boldsymbol{\varphi} = \mathbf{P}_{ss}^{-1} \mathbf{P}_{s\psi} \boldsymbol{\psi} \, ,
  \label{eq:flux_to_source}
\end{equation}
where  $\mathbf{P}_{ss}$ is diagonal with elements $\braket{p^s_j|p^s_j}$ and $\mathbf{P}_{s\psi}$ has elements $\braket{p^s_i| p^{\psi}_j}$. By symmetry, one can also write
\begin{equation}
  \boldsymbol{\psi} = \mathbf{P}_{\psi \psi}^{-1} \mathbf{P}_{\psi s} \boldsymbol{\varphi} \, ,
    \label{eq:source_to_flux}
\end{equation}
and substitution of this expression into Eq.~(\ref{eq:flux_to_source}) yields
\begin{equation}
    \boldsymbol{\varphi} = \mathbf{P}_{ss}^{-1} \mathbf{P}_{s\psi} \mathbf{P}_{\psi \psi}^{-1} \mathbf{P}_{\psi s} \boldsymbol{\varphi} \, .
\end{equation}
Were the bases normalized (so that $\braket{p_i|p_i} = 1$), the matrices $\mathbf{P}_{ss}$ and $\mathbf{P}_{\psi \psi}$ would become the identity, and we would have
\begin{equation}
    \boldsymbol{\varphi} =  \mathbf{P}_{s\psi}  \mathbf{P}_{\psi s} \boldsymbol{\varphi} \, ,
\end{equation}
which is satisfied only if $\mathbf{P}_{s\psi} =  \mathbf{P}_{\psi s}^{-1}$.  This generally is not true, which means that the scalar flux moments (and scattering sources) are not generally preserved.  However, if all bases include the constant function as the first basis function, then averages are preserved in the transformation. (The nonlinear methods referenced earlier are nonlinear because they iteratively adjust the angular flux representation or underlying source representation to preserve the first Legendre moment\cite{walters1996asp, mathews1994ecs, wareing1997exponential, maginot2012nnm}.)

\subsection{An Exponential Basis}

Although polynomial representations often yield simple expressions and efficient numerical implementations, the solutions they represent may not always be well modeled by polynomials.  In particular, linear representations can lead to solutions that are locally negative despite all source terms being strictly positive.

Given an incident flux $\hat{\psi}$ on the left of our cell and a linear source $\tilde{s}(x) = s_0 p_0(x) + s_1 p_1(x)$, the resulting flux is exactly
\begin{equation}
\psi{\left(x \right)} = \left(\hat{\psi} - \frac{s_{0}}{\Sigma} + \frac{s_{1}}{\Sigma} + \frac{2 \mu s_{1}}{\Delta \Sigma^{2}}\right) e^{- \frac{\Sigma x}{\mu}} + \frac{s_{0}}{\Sigma} - \frac{s_{1}}{\Sigma} + \frac{2 s_{1} x}{\Delta \Sigma} - \frac{2 \mu s_{1}}{\Delta \Sigma^{2}} \, .
\label{eq:reference_linear}
\end{equation}
This expression has three distinct functional forms: a constant, a line, and an exponential $e^{-\Sigma x/\mu}$. (This explains why a two-term scheme cannot generally provide an exact solution.)

However, the exact solution does suggest a better, two-term approach: much more of the solution depends on the exponential than the linear term.  In fact, if $s_1 = 0$, then the solution is exact with only a constant and exponential. Hence, we choose for the flux basis the constant function and the exponential function.  However, to make them mutually orthogonal, we end up with
\begin{equation}
  e_0(x) = 1 \, \quad \text{and} \quad
  e_1(x) = \left (1 + \frac{\mu}{\Delta \Sigma} \right) e^{- \frac{\Sigma x}{\mu}} - \frac{\mu}{\Delta \Sigma}  \, .
\end{equation}
The substitution of $\tilde{\psi} = \psi_0 e_0(x) + \psi_1 e_1(x)$ and $\tilde{s}(x) = s_0 p_0(x) + s_1 p_1(x)$ into Eq.~(\ref{eq:projected_general}) with shifted-Legendre test functions leads to
\begin{equation}
\left[
\begin{matrix}
\Delta \Sigma + \mu & \mu e^{- \frac{\Delta \Sigma}{\mu}} - \frac{\mu^{2}}{\Delta \Sigma} + \frac{\mu^{2} e^{- \frac{\Delta \Sigma}{\mu}}}{\Delta \Sigma}\\- \mu & - \mu + \frac{\mu^{2}}{\Delta \Sigma} - \frac{\mu^{2} e^{- \frac{\Delta \Sigma}{\mu}}}{\Delta \Sigma}
\end{matrix}
\right]
\left [ \begin{matrix}
 \psi_0 \\
 \psi_1 \\
\end{matrix} \right ] =
\left [ \begin{matrix}
 \Delta s_0 + \mu \psi_{\text{up}} \\
  \frac{\Delta s_1 }{3} -\mu \psi_{\text{up}} \\
\end{matrix} \right ] \, ,
\label{eq:ex_equations}
\end{equation}
the solution of which is
\begin{equation}
   \psi_{EX} = \left(- \frac{\Delta s_{1}}{3 \mu} + \psi_{up} - \frac{s_{0}}{\Sigma} - \frac{s_{1}}{3 \Sigma}\right) e^{- \frac{\Sigma x}{\mu}} + \frac{s_{0}}{\Sigma} + \frac{s_{1}}{3 \Sigma} \, .
\end{equation}
When $s_1 = 0$, this solution is exact (as expected).  Moreover, it produces a positive outgoing flux as long as $|s_1| < 3s_0$.  For a slab cell, one might wonder what conditions might break this positivity constraint.  Surely, the source would need to fall rapidly from one edge to the next, perhaps exponentially.  A possibly ``reasonable'' bounding case is to assume a source shape equal to the analytic scalar flux in a pure absorber due to a unit isotropic plane source at the left edge, i.e., $\phi_a(x) = E_1(\Sigma x)/2$.  The value of $|s_1|/s_0 = |3\braket{p_1|\phi_a}|/\braket{p_0|\phi_a}$ approaches 3 in the limit of large $\Sigma$ or $\Delta$.  Hence, it seems improbable to observe a negative outgoing flux $\psi_{EX}(\Delta)$.

To determine the resulting scalar-flux moments from Eq.~(\ref{eq:flux_to_source}), recall that $\mathbf{P}_{ss}$ is a diagonal matrix with the $i$th diagonal element equal to $\frac{\Delta}{2i+1}$.  Moreover, by construction $p_0 = e_0$, so $\braket{e_1|p_0} = \braket{p_1|e_0} = 0$, and we need only
\begin{equation}
   \braket{p_1|e_1} = - \frac{\mu}{\Sigma} - \frac{\mu e^{- \frac{\Delta \Sigma}{\mu}}}{\Sigma} + \frac{2 \mu^{2}}{\Delta \Sigma^{2}} - \frac{2 \mu^{2} e^{- \frac{\Delta \Sigma}{\mu}}}{\Delta \Sigma^{2}} \, ,
\end{equation}
to compute
\begin{equation}
    \varphi_0 = \psi_0 \, ,
\end{equation}
and
\begin{equation}
    \varphi_1 = \frac{3\braket{p_1|e_1}}{\Delta} \psi_1 \, .
\end{equation}
This latter expression represents the least-squares fit of the exponential basis function $e_1(x)$ to the shifted-Legendre polynomial $p_1(x)$.

For a strictly linear source, the local errors in the incident, outgoing, and average angular flux of this exponential-discontinuous method are
\begin{equation}
 \displaystyle \epsilon_{in} = \displaystyle - \frac{s_{1} \tau}{3 \Sigma \mu} + O\left(\tau^{5}\right) \, ,
\end{equation}
\begin{equation}
 \displaystyle \epsilon_{out} = - \frac{s_{1} \tau^{3}}{36 \Sigma \mu^{3}} + \frac{s_{1} \tau^{4}}{60 \Sigma \mu^{4}}  + O\left(\tau^{5}\right) \, ,
\end{equation}
and
\begin{equation}
 \displaystyle \epsilon_{avg} =  \frac{s_{1} \tau^{2}}{36 \Sigma \mu^{2}} - \frac{s_{1} \tau^{3}}{60 \Sigma \mu^{3}} + O\left(\tau^{4}\right)  \, ,
\end{equation}
where $\tau = \Sigma\Delta$ is the optical cell width.  The leading error terms for EX are identical to those for LD reported in Eqs. (\ref{eq:ld_error_in}--\ref{eq:ld_error_avg}).  The coefficients for the next highest term are smaller for EX, though, in practice, cancellation of errors may render the difference negligible.

\section{Results and Analysis}

The exponential-discontinuous (EX) method developed above was implemented along the LD, step-characteristic (SC), and linear-characteristic (LC) schemes; here, LC refers to a characteristic solution with a strictly linear source term.  For the sample results to follow, a 16-angle, Gauss-Legendre quadrature was used.  To solve the systems of equations, GMRES was used with a relative tolerance of $10^{-5}$ (with no preconditioning).

To test the EX scheme, a slab representing a deep-penetration application was defined.  The slab is 10 cm thick and surrounded by vacuum.  The first region consists from $x=0$ to $x = 4$ cm contains a uniform, isotropic source of 100 1/s and has $\Sigma_t = 1$ cm$^{-1}$.  The second region is from $x = 4$ to $x = 6$ cm and represents a sourceless shield with $\Sigma_t = 20$ cm$^{-1}$. The final region is also sourceless with $\Sigma_t = 1$ cm$^{-1}$.  Throughout, $\Sigma_s = 0.4 \Sigma_t$.

\begin{figure}[ht] 
  \centering
  \includegraphics[width=0.49\textwidth]{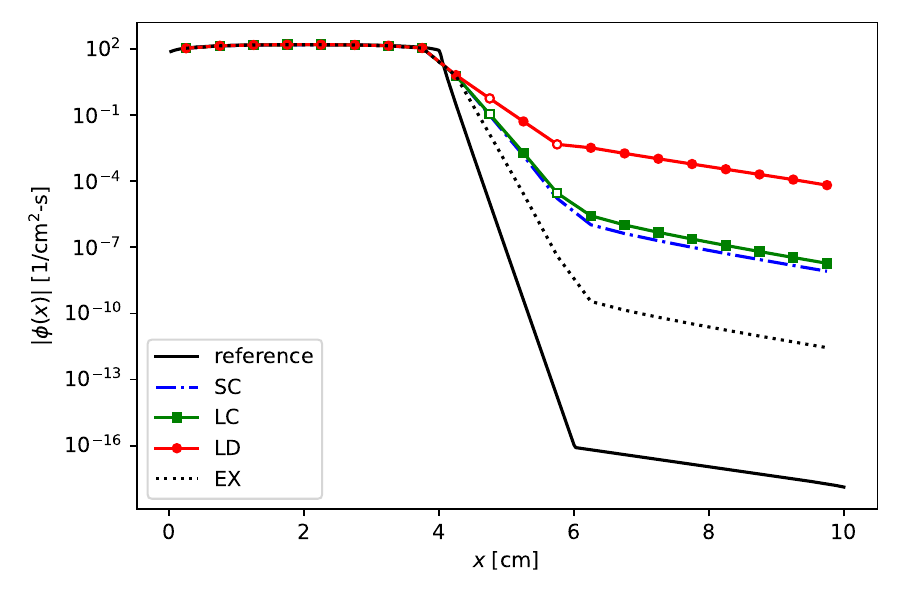}
  \caption{Results for $\Delta = 5$ cm ($\tau = 10$ between 4 and 6 cm).  White markers indicate a negative flux value.}
  \label{fig:results_10}
\end{figure}

\begin{figure}[ht] 
  \centering
  \includegraphics[width=0.49\textwidth]{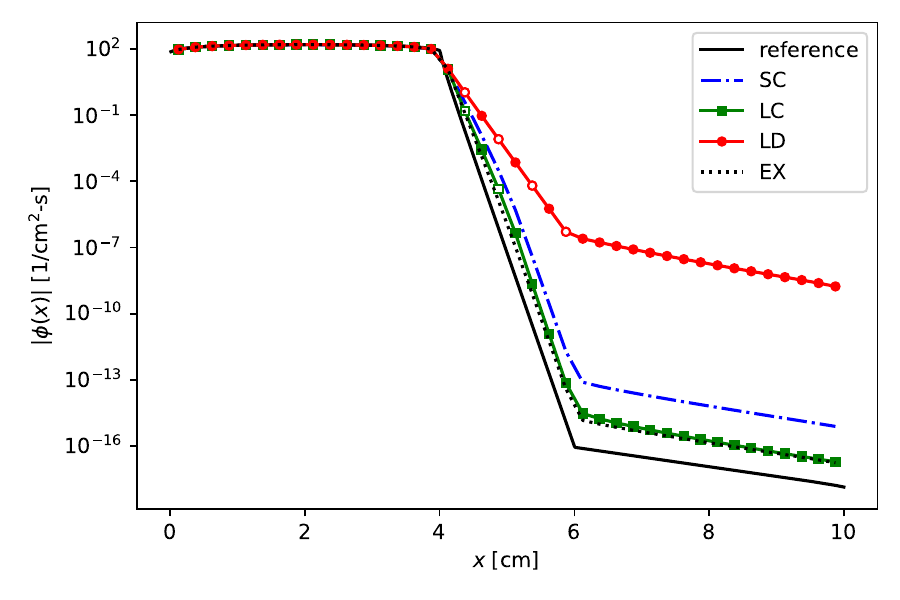}
  \caption{Results for $\Delta = 2.5$ cm ($\tau = 5$ between 4 and 6 cm). White markers indicate a negative flux value.}
  \label{fig:results_5}
\end{figure}

Figures~\ref{fig:results_10} and \ref{fig:results_5} show the scalar flux for a fixed $\Delta = 5$ cm and $\Delta = 2.5$ cm, respectively; with these $\Delta$ values, the optical width of the shield cells is 10 and 5, respectively.  A fine-mesh, LD solution was used as a reference.  For both cases, the EX scheme exhibits the best performance, although for smaller $\Delta$, LC will overtake EX slightly.  Of course, the errors for EX in Figure~\ref{fig:results_10} are still several orders of magnitude---much too large if ``converged'' accuracy is required.  However, for applications in which a ``good enough'' solution is somewhere between the SC (or step-difference) and fully converged, then EX may represent a useful alternative.

That the scalar flux diminishes so rapidly through the shield makes the use of a linear source approximation dubious at best. The equivalent nonlinear techniques described in Refs. \cite{mathews1994ecs} and \cite{walters1996asp} fix this issue by replacing the linear source with one of the form $a e^{bx}$ having the same moments.  It may be possible to define an alternative exponential function {\it a priori} (or to extend the basis) so that sharp attenuation is better captured.

\section{Conclusions}

A strictly linear, exponential-discontinuous EX scheme was presented that exhibits (1) positive angular flux solutions as long the source exhibits a dropoff no worse than $E_1(\Sigma x)$, (2) reasonably good accuracy for optically-thick cells, (3) and local errors of the same order but with potentially smaller coefficients than LD.  Although life is easy in slab geometry, the basic idea of EX should extend naturally to 2-D and 3-D Cartesian grids, though the nature of solutions in higher dimensions is expected to impact the local accuracy and positivity of EX.  This extension and an analysis of additional bases are the subject of ongoing work.

\section{Acknowledgments}

The author thanks E. Sosnovksy for helpful discussions and for suggesting the sample problem.

\bibliographystyle{ans}
\bibliography{bibliography}
\end{document}